\documentclass[reqno,12pt]{amsart}
\usepackage{bbm}
\usepackage{amsmath}
\usepackage{txfonts}
\usepackage{stmaryrd}
\usepackage{amssymb}
\usepackage{amsfonts}
\usepackage{times}
\usepackage{mathrsfs}
\usepackage{amsthm}
\usepackage{enumerate}
\usepackage{framed}
\usepackage{lipsum}
\usepackage{color}

\def\B{\mathcal{B}}
\def\D{\mathbb{D}}

\def\C{\mathbb{C}}

\def\N{\mathbb N}

\def\f{\frac}

\def\B{\mathbb{B}}

\def\a{\alpha}

\def\msk{\medskip}

\def\bege{\begin{equation}} \def\ende{\end{equation}}
  \def\a{\alpha} 
\def\b{\beta}   
 
\def\begr{\begin{eqnarray}} \def\endr{\end{eqnarray}}

\def\bege{\begin{equation}} \def\ende{\end{equation}}
\def\begr{\begin{eqnarray}} \def\endr{\end{eqnarray}}
\def\bnum{\begin{enumerate}} \def\enum{\end{enumerate}}

\begin{document}

\title[Difference of composition operators]{Difference of composition operators on weight Bergman spaces  with doubling weight}
 \author{Yecheng Shi and  Songxiao Li$^*$}

\address{Yecheng Shi\\    School of Mathematics and Statistics, Lingnan Normal University,
     Zhanjiang 524048, Guangdong, P. R. China}\email{ 09ycshi@sina.cn}

\address{Songxiao Li\\ Institute of Fundamental and Frontier Sciences, University of Electronic Science and Technology of China,
610054, Chengdu, Sichuan, P. R. China.   } \email{jyulsx@163.com}

\begin{abstract}    In this paper, some characterizations for the  compact difference of composition operators on Bergman spaces $A^p_\omega$ with doubling weight are given, which extend Moorhouse's characterization for the difference of composition operators on the weighted Bergman space $A^2_\a$.
\thanks{*Corresponding author.}
\thanks{This  project was partially supported by  NNSF of China (No. 11901271 and No. 11720101003) and a grant of Lingnan Normal University (No. 1170919634)}
\end{abstract}

\maketitle

\section{Introduction}
Let $\mathbb{\D}$ be the  the unit disc and $H(\mathbb{\D})$ be the class of analytic functions on $\D$.
 Let $\varphi$ be an analytic self-map of $\D$. The map $\varphi$ induces a composition operator $C_\varphi$ on $H(\D)$, which is defined by
$C_\varphi f=f\circ\varphi$. We refer to \cite{cm,Sh1} for various aspects on the theory of composition operators acting on analytic function spaces.

A function $\omega:\D\to[0,\infty)$, integrable over $\D$, is called a weight. It is radial if $\omega(z)=\omega(|z|)$ for all $z\in\D$. For $0<p<\infty$ and a radial weight $\omega$, the weighted Bergman space $A^p_\omega$ is the space of all $f\in H(\D)$ such that
$$\|f\|_{A^p_\omega}^p=\int_{\D}|f(z)|^p\omega(z)dA(z)<\infty,$$
where $dA(z)$ is the normalized Lebesgue are measure on $\D$. As
usual, $A^p_\a$ stands for the classical weighted Bergman space induced by the standard radial weight $\omega(z)=(1-|z|^2)^\a$, where $-1<\a<\infty$.
$A^p_\omega$ equipped with the norm $\|\cdot\|_{A^p_\omega}$ is a Banach space for $1\leq p<\infty$ and a complete metric space for $0<p<1$ with respect to the translation-invariant metric $(f,g)\mapsto\|f-g\|_{A^p_\omega}$.

For a radial weight $\omega$, we assume throughout the paper that $\widehat{\omega}(r)=\int_{r}^1\omega(s)ds$ for all $0\leq r<1$. A radial weight $\omega$  belongs to $\widehat{\mathcal{D}}$ if there exists a constant $C=C(\omega)>1$ such that
$$\widehat{\omega}(r)\leq C\widehat{\omega}({\frac{1+r}{2}})$$
for all $0\leq r<1$. If there exist $K=K(\omega)>1$ and $C=C(\omega)>1$ such that
$$\widehat{\omega}(r)\geq C\widehat{\omega}(1-\frac{1-r}{K}),~~~\mbox{~~}0\leq r<1,$$
then we say that $\omega\in\check{\mathcal{D}}.$ We write $\mathcal{D}=\widehat{\mathcal{D}}\cap\check{\mathcal{D}}$. For some properties of these classes of weights, see \cite{P,PR1,PR2,PR3,PR4,PR5,PRS} and the references therein.

Efforts to understand the topological structure of the space of composition operators in the operator norm topology have led to the study of the difference operator $C_\varphi-C_\psi$~of two composition operators induced by analytic self-maps $\varphi,\psi$ of $\D$.
By Littlewood's subordination principle, all
composition operators, and hence all differences of two composition operators, are
bounded on all Hardy space $H^p$ and   weighted Bergman spaces $A^p_\a$. Thus
the question of when the operator $C_\varphi-C_\psi$ is compact naturally arises. Shapiro
and Sundberg \cite{SS} raised and studied such a question on   Hardy spaces, motivated by the isolation phenomenon observed by Berkson \cite{Be}. After that, such related
problems have been studied between several spaces of analytic functions by many
authors. See, for example, \cite{G,NS,SL} on Hardy spaces and \cite{CKP2,CKP1,KW,LS,Mo,S1,S2,SLD}
on weighted Bergman spaces.

In 2005, Moorhouse \cite{Mo} characterized the compact difference of composition operators on  weighted Bergman spaces $A^2_\a$ by angular derivative cancellation property. More precisely, she showed that $C_\varphi-C_\psi$ is compact on $A^2_\a$ if and only if
\begin{equation}
\lim_{|z|\to1}\left(\frac{1-|z|^2}{1-|\varphi(z)|^2}+\frac{1-|z|^2}{1-|\psi(z)|^2}\right)\rho(\varphi(z),\psi(z))=0.\label{mo}
\end{equation}
We remark here that this characterization has been extended not only to higher dimensional balls and polydisks, but also to general parameter $p$, see \cite{CKP2,CKP1,LS}.

It is known that all composition operator and hence all differences of two composition operators, are bounded on $A^p_\omega$ for $\omega\in\widehat{\mathcal{D}}$~(see \cite{PR3}).
In this paper we   extend Moorhouse's characterization as well as some related properties to  weighted Bergman spaces $A^p_\omega$, whenever $\omega\in\mathcal{D}$. The approach employed in the proof of the main results of this paper follows the guideline of \cite{CKP1,KM, Mo}, however a good number of steps cannot adapted straightforwardly and need substantial modifications.

Our main  result (Theorem 12)~is a characterization for compact combination of two composition operators. As a consequence we obtain that the Moorhouse's characterization for compact difference (\ref{mo}) remains valid when $0<p<\infty$~and ~$\omega\in\mathcal{D}$.
According to this result, the compactness of $C_\varphi-C_\psi:A^p_\omega\to A^p_\omega$ depends neither on $p$ nor $\omega$, whenever $0<p<\infty$~and $\omega\in\mathcal{D}$. The key ingredient for obtaining the previously mentioned results is the characterization of the $p$-Carleson measure  for $A^p_\omega$.

The present paper is organized as follows. In Section 2, we give some notations and preliminary results which will be used later.
In Sections 3, we devote to the question of when a given finite linear combination of
composition operators is compact.     Section 4 is devoted to show that the Moorhouse's characterization for compact difference remains valid when $0<p<\infty$~and ~$\omega\in\mathcal{D}$. We also obtain a characterization for a composition operator to be equal modulo compact operators to a linear combination of composition operators (see Theorem 14).

For two quantities $A$ and $B$, we use the abbreviation
$A\lesssim B$ whenever there is a positive constant $C$ (independent of the associated variables) such that $A\leq CB$.
We write $A\asymp B$, if $A\lesssim B\lesssim A$.

\section{prerequisites}\vspace{0.2truecm}

In this section we provide some basic tools for the proofs of the main results in this paper.

\subsection{Pseudo-hyperbolic distance}

We denote by $\sigma_z(w)$ the M\"{o}bius transformation on $\D$ that interchanges the points $0$ and $z$. More explicitly,
$$\sigma_z(w)=\f{z-w}{1-\overline{w}z}.$$
It is well known that $\sigma_z$ satisfies the following properties:
$\sigma_z\circ\sigma_z(w)=w$, and
$$1-|\sigma_z(w)|^2=\frac{(1-|z|^2)(1-|w|^2)}{|1-\overline{w}z|^2}, ~~~~z,w\in\D.$$
For $z,w\in\D$, the pseudo-hyperbolic distance between $z$ and $w$ is defined by
$$\rho(z,w)=|\sigma_z(w)|.$$
 It is also well known that the pseudo-hyperbolic metric have the following strong form of triangle inequality (see \cite{DW}):
 $$\rho(z,w)\leq\f{\rho(z,a)+\rho(a,w)}{1+\rho(z,a)\rho(a,w)}$$
for all $a,z,w\in\D$.  For $z\in \D$ and $r>0$, the pseudo-hyperbolic disk at $z\in\D$ with radius $r\in(0,1)$ is given by
 $$\triangle(z,r)=\{w\in\D:\rho(z,w)<r\}.$$
 Note that $\triangle(z,r)$ is Euclidean disk with center and radius given by
$$ c=\f{(1-r^2)z}{1-r^2|z|^2}~~,\mbox{~~~~}~~t=\f{1-|z|^2}{1-r^2|z|^2}r.$$
For $w\in\triangle(z,r)$, it is geometrically clear that
$$|c|-t\leq |w|\leq |c|+t.$$
Therefore,
$$\frac{(1-|z|)(1-r|z|)(1-r)}{1-r^2|z|^2}\leq 1-|w|\leq\frac{(1-|z|)(1+r|z|)(1+r)}{1-r^2|z|^2},$$
and $|w|\to1$ uniformly in $w\in\triangle(z,r)$,~as~$|z|\to1.$

\subsection{Basic properties of weights}

The following two lemmas contains basic properties of weights in the class $\widehat{\mathcal{D}}$ and $\check{\mathcal{D}}$ and will be frequently used in the sequel. For a proof of the first lemma, see \cite[Lemma 2]{P}. The second one can be proved by similar arguments.\msk

\noindent{\bf Lemma A.}  {\it Let $\omega$ be a radial weight. Then the following statements are equivalent:

(i)~$\omega\in\widehat{\mathcal{D}}$;

(ii)~There exist $C=C(\omega)>0$ and $\b=\b(\omega)>0$ such that
 $$\widehat{\omega}(r)\leq C\left(\frac{1-r}{1-t}\right)^{\b}\widehat{\omega}(t),~~0\leq r\leq t<1;$$

 (iii)~There exist $\gamma=\gamma(\omega)>0$ such that
 $$\int_{\D}\frac{dA(z)}{|1-\overline{\zeta}z|^{\gamma+1}}\asymp \frac{\widehat{\omega}(\zeta)}{(1-|\zeta|)^\gamma},~~\zeta\in\D.$$
}\msk

\noindent{\bf Lemma B.}  {\it Let $\omega$ be a radial weight. Then $\omega\in\check{\mathcal{D}}$
if and only if
there exist $C=C(\omega)>0$ and $\alpha=\alpha(\omega)>0$ such that
 $$\widehat{\omega}(t)\leq C\left(\frac{1-t}{1-r}\right)^{\alpha}\widehat{\omega}(r),~~0\leq r\leq t<1.$$
}\msk

The following equivalent norm will be used in our proof, see \cite[Lemma 5]{PR5}.\msk

\noindent{\bf Lemma C.}  {\it Let $0<p<\infty$, $\omega\in\mathcal{D}$ and $-\alpha<\gamma<\infty$, where $\alpha=\alpha(\omega)>0$ is that of Lemma B. Then
 $$\int_{\D}|f(z)|^p(1-|z|^2)^\gamma\omega(z)dA(z)\asymp\int_{\D}|f(z)|^p(1-|z|^2)^{\gamma-1}\widehat{\omega}(z)dA(z),~~f\in H(\D).$$
}\msk

The following estimate plays an important role in this paper and will be frequently used in the sequel.\msk

\noindent{\bf Lemma 1.}  {\it Let $\varphi$ be an analytic self-map of~$\D$ and $\omega\in\mathcal{D}$. Then
\begr
\left(\frac{1-|z|}{1-|\varphi(z)|}\right)^{\b+1}\lesssim\frac{\omega(S(z))}{\omega(S(\varphi(z)))}\lesssim\left(\frac{1-|z|}{1-|\varphi(z)|}\right)^{\a+1},\nonumber
\endr
where $\alpha=\alpha(\omega)$ and $\b=\b(\omega)$ are that of Lemma B and Lemma A, respectively.
}\msk

 {\it Proof.}  An application of Lemma A shows that
$$\omega(S(z))\asymp\widehat{\omega}(z)(1-|z|)~~\mbox{~~and~~}~~\omega(S(\varphi(z)))\asymp\widehat{\omega}(\varphi(z))(1-|\varphi(z)|).$$
By Schwarz's Lemma, we have
$$|\varphi(z)|\leq \frac{c-1}{c}+\frac{|z|}{c},~~~\mbox{~~where~~}c=\frac{1+|\varphi(0)|}{1-|\varphi(0)|}.$$
By Lemmas A and B, we get
\begr
\frac{\widehat\omega(z)}{\widehat\omega(\varphi(z))}&=&\frac{\widehat\omega(z)}{\widehat\omega(\frac{c-1}{c}+\frac{|z|}{c})}\cdot\frac{\widehat\omega(\frac{c-1}{c}+\frac{|z|}{c})}{\widehat\omega(\varphi(z))}\nonumber\\
&\gtrsim&\left(\frac{1-|z|}{1-(\frac{c-1}{c}+\frac{|z|}{c})}\right)^\alpha\left(\frac{1-(\frac{c-1}{c}+\frac{|z|}{c})}{1-|\varphi(z)|}\right)^\b\nonumber\\
&\asymp&\left(\frac{1-|z|}{1-|\varphi(z)|}\right)^\b \nonumber
\endr
and
\begr
\frac{\widehat\omega(z)}{\widehat\omega(\varphi(z))}&=&\frac{\widehat\omega(z)}{\widehat\omega(\frac{c-1}{c}+\frac{|z|}{c})}\cdot\frac{\widehat\omega(\frac{c-1}{c}+\frac{|z|}{c})}{\widehat\omega(\varphi(z))}\nonumber\\
&\lesssim&\left(\frac{1-|z|}{1-(\frac{c-1}{c}+\frac{|z|}{c})}\right)^\beta\left(\frac{1-(\frac{c-1}{c}+\frac{|z|}{c})}{1-|\varphi(z)|}\right)^\alpha\nonumber\\
&\asymp&\left(\frac{1-|z|}{1-|\varphi(z)|}\right)^\alpha .   \nonumber
\endr
  The proof is complete.   $\Box$\msk

\noindent{\bf Lemma 2.}  {\it Let $\omega\in\mathcal{D}$. If $0<\lambda<\a(\omega)$, then  $\omega_{\lambda}(\cdot):=\frac{\omega(\cdot)}{(1-|\cdot|)^\lambda}\in\mathcal{D}$ and
$$\widehat{\omega_\lambda}(z)\asymp\frac{\widehat{\omega}(z)}{(1-|z|)^\lambda},   ~~\mbox{~~~for all~~}~~z\in\D.$$
}

 {\it Proof.} An integration by parts shows that
\begr
\widehat{\omega_\lambda}(r)
=\frac{\widehat{\omega}(r)}{(1-r)^\lambda}+\lambda\int_r^1\widehat\omega(t)(1-t)^{-1-\lambda}dt.   \nonumber
\endr
Therefore, by Lemmas A and B, we have
\begr
\widehat{\omega_\lambda}(r)\gtrsim\frac{\widehat{\omega}(r)}{(1-r)^\lambda}+\lambda\frac{\widehat{\omega}(r)}{(1-r)^\b}\int_r^1(1-t)^{\beta-1-\lambda}dt\gtrsim\frac{\widehat{\omega}(r)}{(1-r)^\lambda}\nonumber
\endr
and
\begr
\widehat{\omega_\lambda}(r)\lesssim\frac{\widehat{\omega}(r)}{(1-r)^\lambda}+\lambda\frac{\widehat{\omega}(r)}{(1-r)^\a}\int_r^1(1-t)^{\a-1-\lambda}dt\lesssim\frac{\widehat{\omega}(r)}{(1-r)^\lambda}.\nonumber
\endr
Thus,
$$\widehat{\omega_\lambda}(z)\asymp\frac{\widehat{\omega}(z)}{(1-|z|)^\lambda}~~\mbox{~~~for all~~}~~z\in\D.$$
By Lemmas A and B, $\omega_{\lambda}\in\mathcal{D}$.   $\Box$

\subsection{Local estimates and test functions}\vspace{0.3truecm}

The following lemmas are crucial in our work and will be   used in this paper.\msk

\noindent{\bf Lemma 3.}  {\it Let $0<p<\infty$, $\omega\in\widehat{\mathcal{D}}$ and $0<r_1<1$ be arbitrary. Denote $\widetilde\omega(\cdot)=\frac{\widehat\omega(\cdot)}{1-|\cdot|}$. Then there exists $0<r_2<1$ and a constant $C=C(p,\omega,r_1,r_2)$ such that
\begr
|f(z)-f(a)|^p\leq C\rho(z,a)^p\frac{\int_{\bigtriangleup(z,r_2)}|f(\zeta)|^p\widetilde\omega(\zeta)dA(\zeta)}{\omega(S(z))}\nonumber
\endr
for all $a\in\D$, $z\in \triangle(a,r_1)$ and $f\in A^p_\omega$.
}\msk

 {\it Proof.} Let $a\in\D$, $0<r_1<1$, $r:=\frac{2r_1}{1+r_1^2}$, $\delta:=\frac{2r}{1+r^2}$,$r_2:=\frac{2\delta}{1+\delta^2}$  and $z\in \triangle(a,r_1)$ be fixed. Consider $g_a:=f\circ\sigma_a$.
Then,
\begr
\big|f(z)-f(a)\big|^p&=&|g_a(\sigma_a(z))-g_a(0)|^p\nonumber\\
&=&|g_a^{\prime}(\eta)|^p|\sigma_a(z)|^p\nonumber\\
&=&|\sigma_a(z)|^p\big|\frac{1}{2\pi}\int_{|\xi|=r}\frac{g_a(\xi)}{(\xi-\eta)^2}d\xi\big|^p\nonumber
\endr
for some $\eta$ with $|\eta|\leq|\sigma_a(z)|<r_1$.
Since $|\xi|=\rho(\sigma_a(\xi),a)=r$, we get $u:=\sigma_a(\xi)\in\triangle(a,\delta)$.
Thus
\begr
\big|f(z)-f(a)\big|^p&\lesssim&\rho(a,z)^p\left(\frac{1}{2\pi}\int_{|\xi|=r}\left|\frac{g_a(\xi)}{(\xi-\eta)^2}\right|d\xi\right)^p\nonumber\\
&\lesssim&\rho(a,z)^p\sup_{u\in\triangle(a,\delta)}|f(u)|^p.   \label{one}
\endr
 Using the subharmonicity of $|f(u)|^p$, $1-|u|\asymp1-|\zeta|$ for $\zeta\in\triangle(u,\delta)$, and
 \begin{equation}
\omega(S(\zeta))\asymp\widehat\omega(\zeta)(1-|\zeta|)\asymp\widetilde\omega(\zeta)(1-|\zeta|)^2, \nonumber
\end{equation}
 we get
\begr
|f(u)|^p&\lesssim& \frac{1}{(1-|u|^2)^2}\int_{\triangle(u,\delta)}|f(\zeta)|^pdA(\zeta)\nonumber\\
&\lesssim& \int_{\triangle(u,\delta)}|f(\zeta)|^p\frac{\widetilde{\omega}(\zeta)}{\omega{(S(\zeta))}}dA(\zeta)\nonumber\\
&\lesssim& \frac{1}{\omega{(S(a))}}\int_{\triangle(a,r_2)}|f(\zeta)|^p\widetilde{\omega}(\zeta)dA(\zeta),
\label{two}
\endr
where we use the fact that $\triangle(u,\delta)\subset\triangle(a,r_2)$ and
 \begin{equation}
\omega(S(a))\asymp \omega(S(\zeta)),\label{ie}
\end{equation}
for $\zeta\in\triangle(a,r_2)$.
Combining (\ref{one}) and (\ref{two}), we obtain
\begr
\big|f(z)-f(a)\big|^p
&\leq& C \rho(z,a)^p\frac{\int_{\triangle(a,r_2)}|f(\zeta)|^p\widetilde{\omega}(\zeta)dA(\zeta)}{\omega(S(a))}.\nonumber
\endr
The proof is complete. $\Box$\msk

By \cite[Lemma 4.30]{Zhu}, for all $a,z,w\in\mathbb D$ with $\rho(z,w)<r$ and any real $s$, we have
$$\left|1-\left(\frac{1-\overline{a}z}{1-\overline{a}w}\right)^s\right|\leq C(s,r)\rho(z,w),$$
and therefore, for all $w,z,a\in\mathbb D$ with $z\in\triangle(a,r)$ and any $s>0$,
$$\left|\frac{1}{(1-\overline{a}z)^s}-\frac{1}{(1-\overline{a}w)^s}\right|\leq C(s,r)\rho(z,w)\left|\frac{1}{(1-\overline{a}z)^s}\right|.
 $$
Although the converse inequality does not hold,   we have the following partial converse inequality (see \cite[Theorem 2.8]{KW} or \cite[Lemma 2.3]{SLD}),
 which is crucial in the proof of the necessary part
  of  Theorems 12 and 14.\msk

\noindent{\bf Lemma D.}  {\it Suppose $s>1$ and $0<r_0<1$. Then there are $N=N(r_0)>1$ and $C=C(s,r_0)$ such that
\begr
&&\left|\frac{1}{(1-\overline{a}z)^s}-\f{1}{(1-\overline{a}w)^s}\right|+
\left|\frac{1}{(1-t_N\overline{a}z)^s}-\f{1}{(1-t_N\overline{a}w)^s}\right|\nonumber\\
&\geq& C\rho(z,w)\left|\frac{1}{(1-\overline{a}z)^s}\right|, \nonumber
\endr
for all $z\in\triangle(a,r_0)$ with $1-|a|<\f{1}{2N}$, $t_N=1-N(1-|a|)$ and $w\in\D$.
}

\subsection{Carleson measure}\vspace{0.3truecm} Let $\mu$ be a finite positive Borel measure on $\D$. $\mu$ is called a $q$-Carleson measure for $A^p_\omega$ if the
identity operator $I_d:A^p_\omega\to L^q(d\mu)$ is bounded,
 i.e. there is a positive constant $C>0$ such that $$\int_{\D}|f(z)|^qd\mu(z)\leq C\|f\|_{A^p_\omega}^q$$ for any $f\in A^p_\omega$.
Also, $\mu$ is called a vanishing $q$-Carleson measure if the
identity operator $I_d:A^p_\omega\to L^q(d\mu)$ is compact.

The characterization  of (vanishing) $q$-Carleson measure for $A^p_\omega$ has been solved for $\omega\in\widehat{\mathcal{D}}$ \cite{PR1,PRS}. It is worth mentioning that   the pseudohyper-bolic disk is not the right one to describe the Carleson measure for $A^p_\omega$ when $\omega\in\widehat{\mathcal{D}}$, since for a fixed $r>0$, the quantity $\omega(\triangle(a,r))$~may equal to zero for some $a$~close to the boundary if $\omega\in\widehat{\mathcal{D}}$ (see \cite{PR2}).
~However, if $\omega\in\mathcal{D}$,  we have the following characterization. The proof is similar with the proof of Theorem 2.1 in \cite{PR1}.  We give the proof here for completeness. \msk

\noindent{\bf Theorem 4.}  {\it   Let $\mu$ be a positive Borel measure on $\D$, $0<p<\infty$, $\omega\in\mathcal{D}$ and $0<r<1$. Then the following assertions hold:

(i) $\mu$ is a $p$-Carleson measure for $A^p_\omega$ if and only if
\begin{equation}\sup_{a\in\D}\frac{\mu(\triangle(a,r))}{\omega(S(a))}<\infty.\label{eb}
\end{equation}

(ii) $\mu$ is a vanishing $p$-Carleson measure for $A^p_\omega$ if and only if
\begin{equation}
\lim_{|a|\to1}\frac{\mu(\triangle(a,r))}{\omega(S(a))}=0.\label{ec}
\end{equation}}\msk

\noindent{\bf Remark.}  { In the above, $\omega(S(a))$~ can be replaced by $\omega(\triangle(a,r)))$ for
any fixed $r\in(0, 1)$~large enough.}\msk

 {\it Proof.} $(i)$ Assume first that $\mu$ is a $p$-Carleson measure for $A^p_\omega$. Consider the test functions
$$f_{a}(z)=\left(\frac{1-|a|^2}{1-\overline{a}z}\right)^{\frac{\gamma+1}{p}},$$
where $\gamma=\gamma(\omega)>0$ is chosen large enough.
Then the assumption together with Lemma A yield
$$\mu(\triangle(z,r))\lesssim \int_{\triangle(z,r)}|f_a(z)|^pd\mu(z)\lesssim
\|f_a\|_{A^p_\omega}^p\lesssim\omega(S(z)).$$

Conversely, assume that   (\ref{eb}) holds. By Fubini's Theorem, Lemma C and the following well known estimate
$$|f(z)|^p\lesssim \frac{1}{(1-|z|^2)^2}\int_{\triangle(z,r)}|f(\zeta)|^pdA(\zeta),~~~z\in\D, $$
 we have
\begr
\int_{\D}|f(z)|^pd\mu(z)&\lesssim&\int_{\D}\left(\frac{1}{(1-|z|^2)^2}\int_{\triangle(z,r)}|f(\zeta)|^pdA(\zeta)\right)d\mu(z)\nonumber\\
&=&\int_{\D}|f(\zeta)|^p\frac{\mu(\triangle(\zeta,r))}{(1-|\zeta|^2)^2}dA(\zeta)\nonumber\\
&\lesssim&\int_{\D}|f(\zeta)|^p\frac{\omega(S(\zeta))}{(1-|\zeta|^2)^2}dA(\zeta)\nonumber\\
&\lesssim&\int_{\D}|f(\zeta)|^p\frac{\widehat{\omega}(\zeta)}{(1-|\zeta|^2)}dA(\zeta)\nonumber\\
&\lesssim&\|f\|_{A^p_\omega}^p.\nonumber
\endr

$(ii)$  Assume first that $\mu$ is a vanishing $p$-Carleson measure for $A^p_\omega$. Following the proof of \cite[Theorem 2.1($ii$)]{PR1}, with Lemma B
in hand, we get
$$\lim_{|a|\to1}\frac{\mu(\triangle(a,r))}{\omega(S(a))}=0.$$

Conversely, assume that   (\ref{ec}) holds.  Denote $\D_s=\{z\in\D:|z|<s\}$ and set
$$d\mu_s(z)=\chi_{s\leq |z|<1}(z)d\mu(z).$$
We claim that $(i)$ implies
$$\|h\|_{L^q_\omega}\leq K_{\mu_s}\|h\|_{A^p_\omega},~~h\in A^p_\omega,$$
where
$$K_{\mu_s}=\sup_{a\in\D}\frac{\mu_{s}(\triangle(a,r))}{\omega(S(a))}.$$
Following the proof of \cite[Theorem 2.1($ii$)]{PR1}, it remains to show that
$$\lim_{s\to1^{-}}K_{\mu_s}=\lim_{s\to1^{-}}\left(\sup_{a\in\D}\frac{\mu_{s}(\triangle(a,r))}{\omega(S(a))}\right)=0.$$

Let $t_r(s)=\frac{s-r}{1-sr}$. After an easy calculation, we get that
$\triangle(a,r)\cap(\D\backslash\D_s)\neq\emptyset$ if and only if $|a|\geq t_r(s).$
It is easy to see that $t_r(s)$ is continuous and increasing on $[r,1)$, and $\lim_{s\to1}t_r(s)=1$.
Thus,
\begr
0=\limsup_{|a|\to1}\frac{\mu(\triangle(a,r)))}{\omega(S(a)}&=&\lim_{s\to1}\sup_{|a|\geq t_r(s)}\frac{\mu(\triangle(a,r)))}{\omega(S(a)}\nonumber\\
&\geq&\lim_{s\to1}\sup_{|a|\geq t_r(s)}\frac{\mu(\triangle(a,r))\cap(\D\backslash\D_s))}{\omega(S(a)}\nonumber\\
&=&\lim_{s\to1}\sup_{a\in\D}\frac{\mu_s(\triangle(a,r))}{\omega(S(a)}.\nonumber
\endr
The proof is complete.    $\Box$\msk

The connection between composition operators and Carleson measures comes from the standard identity
$$\int_{\D}(f\circ\varphi)(z) \omega(z)dA(z)=\int_{\D}f(z)d\nu(z),$$
where $\nu$ denotes the pullback measure defined by
$\nu(E)=\int_{\varphi^{-1}(E)}\omega(z)dA(z),$
for all Borel sets $E\subset\D$. On can easily see from the above  equality that $C_\varphi:A^p_\omega\to A^p_\omega$ is bounded (compact) on $A^p_\omega$ if and if $\nu$ is a (vanishing $p$-Carleson measure) $p$-Carleson measure for $A^p_\omega$.\msk

The following result plays a fundamental role in this study. It can be proved by employing the method used by Moorhouse \cite{Mo}.
\msk

\noindent{\bf Lemma 5.}  {\it Let $\varphi$ be an analytic self-map of $\D$, $\omega\in\mathcal{D}$, and $u$ to be a non-negative, bounded, measurable function on $\D$. Define the measure
$\nu(E)=\int_{E}u(z)\omega(z)dA(z)$ on all Borel subset $E$  of $\D$. If
$$\lim_{|z|\to1}u(z)\frac{1-|z|}{1-|\varphi(z)|}=0,$$
then $\nu\circ\varphi^{-1}$ is a vanishing $p$-Carleson measure for $A^p_\omega$ and hence the inclusion map $I_{p,\omega}:A^p_\omega\to L^p(\nu\circ\varphi^{-1})$  is compact.
}\msk

 {\it Proof.} Fix $r\in (0,1)$. For $a\in\D$, let
$$\epsilon:=\epsilon(a)=\sup_{z\in\varphi^{-1}(\triangle(a,r))}u(z)\frac{1-|z|}{1-|\varphi(z)|}.$$
Using the Schwarz-Pick Theorem, one has
$$\frac{1-|z|}{1-|\varphi(z)|}\leq \frac{1-|\varphi(0)|}{1-|\varphi(0)|}=C<\infty.$$
So that if $\varphi(z)\in \triangle(a,r)$, then
$$1-|z|\leq C(1-|\varphi(z)|)\leq C\frac{(1-|a|)(1-r|a|)(1+r)}{1-r^2|a|^2}.$$
This implies that $|z|\to1$ uniformly in $z\in\varphi^{-1}(\triangle(a,r))$~ as ~$|a|\to1$.
By hypothesis $\epsilon(a)\to0$ as $|a|\to1$.

Now, fix $0<\lambda<\min\{1,\alpha(\omega)\}$. Taking $M$ to be an upper bound of $u$,
we have
\begr
\nu\circ\varphi^{-1}(\triangle(a,r)))&=&\int_{\varphi^{-1}(\triangle(a,r)))}u(z)\omega(z)dA(z)\nonumber\\
&\lesssim&\int_{\varphi^{-1}(\triangle(a,r)))}\frac{\epsilon^\lambda(1-|\varphi(z)|)^\lambda}{(1-|z|)^\lambda}u(z)^{1-\lambda}\omega(z)dA(z)\nonumber\\
&\lesssim&\epsilon^\lambda M^{1-\lambda}(1-|a|)^\lambda\int_{\varphi^{-1}(\triangle(a,r)))}\frac{\omega(z)}{(1-|z|)^\lambda}dA(z).  \nonumber
\endr
Denote $\omega_\lambda(z)=\frac{\omega(z)}{(1-|z|)^\lambda}$. By Lemma 2, we get
$\omega_\lambda\in\mathcal{D}$. Therefore,   $C_\varphi:A^p_{\omega_\lambda}\to A^p_{\omega_\lambda}$ is bounded, that is
\begr
(1-|a|)^\lambda\int_{\varphi^{-1}(\triangle(a,r)))}\frac{\omega(z)}{(1-|z|)^\lambda}dA(z)
&\leq& (1-|a|)^\lambda\omega_\lambda(\triangle(a,r)))\nonumber\\
&\asymp&\widehat{\omega_{\lambda}}(a)(1-|a|)^{1+\lambda}\nonumber\\
&\asymp& \widehat{\omega}(a)(1-|a|) \nonumber\\
&\asymp& \omega(\triangle(a,r)).\nonumber
\endr
Therefore
$$\frac{\nu\circ\varphi^{-1}(\triangle(a,r)))}{\omega(\triangle(a,r))}\lesssim\epsilon(a)^{1-\lambda}$$
for all $a\in\D$, and hence we conclude that $\nu\circ\varphi^{-1}$ is a vanishing $p$-Carleson measure for $A^p_\omega$.
The proof is complete.   $\Box$

\subsection{Angular Derivative}\vspace{0.3truecm} Let $\varphi$ be an analytic self-map of $\D$. We say that $\varphi$ has an angular derivative, denoted by $\varphi^\prime(\zeta)\in\C$, at $\zeta\in\partial\D$ if $\varphi$ has nontangential limit $\varphi(\zeta)\in\partial\D$  such that
\begr
\angle \lim_{\begin{subarray}{l}
 z\to\zeta
\end{subarray}}\frac{\varphi(z)-\eta}{z-\zeta}=\varphi^\prime(\zeta), \nonumber
\endr
where $\angle\lim$~stands for the nontangential limit.
We denote by $F(\varphi)$~the set of all boundary points at which $\varphi$~has finite angular derivatives.
Note from the Julia-Carath{\'e}odory Theorem (see \cite[Theorem 2.44]{cm}) that
$$F(\varphi)=\bigg\{\zeta\in\partial\D:d_\varphi(\zeta):=\liminf_{z\to\zeta}\frac{1-|\varphi(z)|}{1-|z|}<\infty\bigg\}.$$
For $\zeta\in F(\varphi)$,  we call the vector
$$\mathcal{D}(\varphi,\zeta):=(\varphi(\zeta),d_\varphi(\zeta))\in\partial\D\times\mathbb{R}^+$$
the first-order data of $\varphi$ at $\zeta$.

If $\varphi$ and $\psi$ are two analytic self-maps of the disk with finite
angular derivative at $\D$, we  say that~$\varphi$ and $\psi$  have the same first-order data at $\zeta$
if $\mathcal{D}(\varphi,\zeta)=\mathcal{D}(\psi,\zeta)$.\msk

\section{linear combination of composition operators}

For a linear operator $T:X\to Y$, the essential norm of $T$, denoted by $\|T\|_{e,X\to Y}$,  is defined by
 $$\|T\|_{e,X\to Y}=\inf\{\|T-K\|_{X\to Y}:K~ \mbox{~is compact from~} X \mbox{~to~} Y\}.$$ It is obvious that the  operator $T$ is compact if and only if $\|T\|_{e,X\to Y}=0$.\msk

We have the following lower estimates for the essential norm of a linear combination of composition operators acting on
  Bergman spaces with doubling weight.\msk

\noindent{\bf Lemma 6.}  {\it Let $0<p< \infty$ and $\omega\in\widehat{\mathcal{D}}$. Let
$\varphi_1,...,\varphi_n$ be finitely many analytic self-maps of $\D$. Then there is a constant $C>0$ and $\gamma=\gamma(\omega)$ is sufficiently large such that
\begr
\Big\|\sum_{j=1}^n\lambda_jC_{\varphi_j}\Big\|_{e,A^p_\omega}^p\geq C\limsup_{|a|\to1}\Big\|(\sum_{j=1}^n\lambda_jC_{\varphi_j})f_{a}\Big\|_{A^p_\omega}^p,\nonumber
\endr
where $f_{a}(z)=\left(\frac{1-|a|^2}{1-\overline{a}z}\right)^{\frac{\gamma+1}{p}}\omega(S(a))^{-\f{1}{p}}$. }\msk

 {\it Proof.}  Let $K$ be a compact operator on $A^p_\omega$. Consider the operator on $H(\D)$ defined by
$$K_m(f)(z)=f(\f{m}{m+1}z),~~m\in\N.$$
Denote $R_m=I-K_m$. It is easy to see that $K_m$ is compact on $A^p_\omega$ (see \cite[Theorem 15]{PR3}) and
$$\|K_m\|_{A^p_\omega}\leq1, \,~~~~\,  \|R_m\|_{A^p_\omega}\leq2$$ for any positive integer $m$.
Then we have
\begr
2\Big\|(\sum_{j=1}^n\lambda_jC_{\varphi_j})-K\Big\|_{A^p_\omega} &\geq& \Big\|R_m\circ(\sum_{j=1}^n\lambda_jC_{\varphi_j}-K)\Big\|_{A^p_\omega}\nonumber\\
&\gtrsim&\sup_{a\in\D}\Big\|R_m\circ(\sum_{j=1}^n\lambda_jC_{\varphi_j}-K)(f_{a})\Big\|_{A^p_\omega}.\nonumber
\endr
Since $K$ is compact, we can extract a sequence $\{a_i\}\subset\D$ such that $|a_i|\to1$ and $Kf_{a_i}$ converges to some $f\in A^p_\omega$.
So,
\begr
&&\Big\|R_m\circ(\sum_{j=1}^n\lambda_jC_{\varphi_j}-K)(f_{a_i})\Big\|_{A^p_\omega}^p\nonumber\\
&\gtrsim&\Big\|R_m\circ(\sum_{j=1}^n\lambda_jC_{\varphi_j})(f_{a_i})\Big\|_{A^p_\omega}^p
-\Big\|R_m\circ K(f_{a_i})\Big\|_{A^p_\omega}^p\nonumber\\
&\gtrsim&\Big\|(\sum_{j=1}^n\lambda_jC_{\varphi_j})(f_{a_i})\Big\|_{A^p_\omega}^p-\Big\|K_m\circ(\sum_{j=1}^n\lambda_jC_{\varphi_j})(f_{a_i})\Big\|_{A^p_\omega}^p\nonumber\\
&~~&-\Big\|R_m(K(f_{a_i})-f)\Big\|_{A^p_\omega}^p-\Big\|R_m(f)\Big\|_{A^p_\omega}^p.\label{1}
\endr
Since $K_m$ is compact and $\sum_{j=1}^n\lambda_jC_{\varphi_j}$ is bounded on $A^p_\omega$, we have
$K_m\circ(\sum_{j=1}^n\lambda_jC_{\varphi_j})$ is compact on $A^p_\omega$.
Therefore, letting $i\to\infty$ and then using Fatou's Lemma as $m\to\infty$ in (\ref{1}), we have
$$\Big\|\sum_{j=1}^n\lambda_jC_{\varphi_j}-K\Big\|_{A^p_\omega}\gtrsim
\limsup_{i\to\infty}\Big\|(\sum_{j=1}^n\lambda_jC_{\varphi_j})(f_{a_i})\Big\|_{A^p_\omega}.
$$
Therefore,
\begr
\Big\|\sum_{j=1}^n\lambda_jC_{\varphi_j}\Big\|_{e,A^p_\omega}^p\geq C\limsup_{|a|\to1}\Big\|(\sum_{j=1}^n\lambda_jC_{\varphi_j})f_{a}\Big\|_{A^p_\omega}^p.\nonumber
\endr
The proof is complete.   $\Box$\msk

For $M>1$ and $\zeta\in\partial\D$, we denote by $\Gamma_{M,\zeta}$ the $\zeta$-curve consisting of points $|z-\zeta|=M(1-|z|^2)$, the boundary of a non-tangential approach region with vertex at $\zeta$. We will use the notation $``\lim_{\Gamma_{M,\zeta}}"$ to indicate a limit taken as $z\to\zeta$ along the stardoard leg of~$\Gamma_{M,\zeta}$. The following result taken from \cite{KM}.\msk

\noindent{\bf Lemma E.}  {\it Let $\varphi$ and $\psi$ be analytic self-maps of $\D$. Then the following equality
\begr
\lim_{M\to\infty}\lim_{\begin{subarray}{l}
z\to\zeta\\
z\in\Gamma_{M,\zeta}
\end{subarray}
}\frac{1-|\varphi(z)|^2}{1-\overline{\varphi(z)}\psi(z)}=\left\{
\begin{array}{ll}
1, &{\rm if~}\zeta\in F(\varphi)\mbox{~and~}\mathcal{D}(\varphi,\zeta)=\mathcal{D}(\psi,\zeta)\\
0, &{\rm otherwise}
\end{array}
\right.
\endr
holds for $\zeta\in F(\varphi)$.
}\msk

We are now ready to establish a lower estimate for the essential norm of a general linear combination of composition operators acting on $A^p_\omega$ when  $\omega\in\widehat{\mathcal{D}}$. Let $\varphi_1,...,\varphi_n$ be finitely many analytic self-maps of $\D$. For $\varphi\in F(\varphi_i)$, we denote by $J_{\zeta}(i)$ the set of all indices $j$ for
which $\zeta\in F(\varphi)$ and $\varphi_i$ and $\varphi$ have the same first-order data at $\zeta$.\msk

\noindent{\bf Theorem 8.}  {\it Let $0<p< \infty$ and $\omega\in\widehat{\mathcal{D}}$. Let
$\varphi_1,...,\varphi_n$ be finitely many analytic self-maps of $\D$. Then there is a constant $C(p,\omega)>0$ such that
\begr
\Big\|\sum_{j=1}^n\lambda_jC_{\varphi_j}\Big\|_{e,A^p_\omega}^p\geq C\max_{1\leq i\leq n}\left(\left|\sum_{j\in J_{\zeta}(i)}\lambda_j\right|^p\f{1}{d_{\varphi_i}(\zeta)^{\b+1}}\right)
\endr
for all $\zeta\in \partial\D$ and $\lambda_1,...,\lambda_n\in\C$. In case $\zeta\notin F(\varphi_i)$
the quantity inside the parenthesis above is to be understood as 0.
}\msk

 {\it Proof.} We denote $T:=\sum_{j=1}^n\lambda_jC_{\varphi_j}$ and $f_{a}(z)=\left(\frac{1-|a|^2}{1-\overline{a}z}\right)^{\frac{\gamma+1}{p}}\omega(S(a))^{-\f{1}{p}}$, for $a\in\D$ and $\gamma$~is that of Lemma A. Fix any index $i$ such that $\zeta\in F(\varphi_i)$.  We have $|\varphi_i(z)|\to1$ as
$z\to\zeta$ along any $\Gamma_{M,\zeta}$ which is a restriced $\zeta$-curve. So, by Lemma 6, we obtain
\begr
\|T\|_{e,A^p_\omega}&\gtrsim&\sup_{M}\left(\lim_{\begin{subarray}{l}
z\to\zeta\\
z\in\Gamma_{M,\zeta}
\end{subarray}}
\|Tf_{\varphi_i(z)}\|_{A^p_\omega}^p\right).\nonumber
\endr
Meanwhile, note that
\begr
\|Tf_{\varphi_i(z)}\|_{A^p_\omega}^p&\geq&|Tf_{\varphi_i(z)}(z)|^p\omega(S(z))\nonumber\\
&=&\left|\sum_{j=1}^n\lambda_j\left(\frac{1-|\varphi_{i}(z)|^2}{1-\overline{\varphi_i(z)}\varphi_j(z)}\right)^{\f{\gamma+1}{p}}\right|^p
\frac{\omega(S(z))}{\omega(S(\varphi_i(z)))}.\nonumber
\endr
Thus, applying Lemma E, Lemmas 1 and   6, we get the desired result. $\Box$\msk

By Theorem 8, we immediately yield  the following three corollaries  for the compactness of linear combinations.   \msk

\noindent{\bf Corollary 9.}  {\it Let $0<p< \infty$ and $\omega\in\widehat{\mathcal{D}}$. Let
$\varphi_1,...,\varphi_n$ be finitely many analytic self-maps of $\D$. If $\sum_{j=1}^n\lambda_jC_{\varphi_j}$ is compact on $A^p_\omega$, then
\begr
\sum_{\begin{subarray}{l}
\zeta\in F(\varphi_j)\\
\mathcal{D}(\varphi_j,\zeta)=(\eta,s)
\end{subarray}}\lambda_j=0    \nonumber
\endr
for all $\zeta\in \partial\D$ and $(\zeta,s)\in\partial\D\times\mathbb{R}_+$.
}\msk

\noindent{\bf Corollary 10.}  {\it Let $0<p< \infty$ and $\omega\in\widehat{\mathcal{D}}$. Let
$\varphi,\psi$ be analytic self-maps of $\D$. Suppose both $C_\varphi$ and $C_\psi$ are not compact on $A^p_\omega$. If $aC_{\varphi}+bC_{\psi}$ is compact on $A^p_\omega$, then the following statements hold:

(i)~$a+b=0$;

(ii)~$F(\varphi)=F(\psi)$;

(iii)~$\mathcal{D}(\varphi,\zeta)=\mathcal{D}(\psi,\zeta)$ for each $\zeta\in F(\varphi)$.
}\msk

\noindent{\bf Corollary 11.}  {\it Let $0<p< \infty$ and $\omega\in\widehat{\mathcal{D}}$. Let
$\varphi,\varphi_1,...,\varphi_n$ be finitely many analytic self-maps of $\D$. If $C_{\varphi}-C_{\varphi_1}-C_{\varphi_2}-\cdots-C_{\varphi_n}$ is compact on $A^p_\omega$, then the following statements hold:

(i)~$F(\varphi_1),\cdots,F(\varphi_n)$ are pairwise disjoint and $F(\varphi)=\cup_{j=1}^nF(\varphi_j)$

(ii)~$\mathcal{D}(\varphi,\zeta)=\mathcal{D}(\varphi_j,\zeta)$ at each $\zeta\in F(\varphi_j)$ for $j=1,\cdots,n$.
}\msk

\section{compact difference and further related results}

We have the following characterization for compact linear combinations of two composition operators.\msk

\noindent{\bf Theorem 12.}  {\it Let $0<p< \infty$ and $\omega\in\mathcal{D}$. Suppose
$\varphi$ and $\psi$ be analytic self-maps of $\D$. Then $\lambda_1C_\varphi+\lambda_2C_\psi$ is compact on $A^p_\omega$ if and only if
either one of the following two conditions holds:

(i) Both $C_\varphi$ and $C_\psi$ are compact;

(ii) $\lambda_1+\lambda_2=0$ and
\begin{equation}
\lim_{|z|\to1}\left(\frac{1-|z|^2}{1-|\varphi(z)|^2}+\frac{1-|z|^2}{1-|\psi(z)|^2}\right)\rho(\varphi(z),\psi(z))=0.\label{cc}
\end{equation}
}\msk

 {\it Proof.} Suppose that $\lambda_1C_\varphi+\lambda_2C_\psi$ is compact on $A^p_\omega$. Note that if $(i)$ fails, then one of $C_\varphi$ and $C_\psi$ is not compact on $A^p_\omega$. We may assume that both $C_\varphi$ and $C_\psi$ are not compact on
$A^p_\omega$ and show $(ii)$. By Corollary 10, we have $\lambda_1+\lambda_2=0$ and hence we may assume that $\lambda_1=1$ and $\lambda_2=-1$. We assume that (\ref{cc})  does not hold. Then there exists a sequence $\{z_n\}\subset\D$ with $|z_n|\to1$ such that either
$$a_n=\frac{1-|z_n|}{1-|\varphi(z_n)|}\rho(\varphi(z_n),\psi(z_n))$$
or
$$b_n=\frac{1-|z_n|}{1-|\psi(z_n)|}\rho(\varphi(z_n),\psi(z_n))$$
does not converge to zero.
By passing to a subsequence, we may assume that $\lim_{n\to\infty}a_n=a$ and $\lim_{n\to\infty}b_n=b$ exist and that one is non-zero. Without loss of generality we may further assume that $a\neq0$.
Again by passing to a subsequence, we may assume that $c=\lim_{n\to\infty}|\varphi(z_n)|$ exist. Since $a\neq0$, we have $c=1$.
Thus, we may assume that $|z_n|\to1$, $|\varphi(z_n)|\to1$ and $a=\lim_{n\to\infty}a_n$ exists and non-zero.
For $u\in\D$, consider the test functions
$$g_{u}(z)=\left(\frac{1-|u|^2}{1-\overline{u}z}\right)^{\frac{\gamma+1}{p}}\omega(S(u))^{-\f{1}{p}},$$
and
$$h_{u}(z)=\left(\frac{1-|u|^2}{1-t_N\overline{u}z}\right)^{\frac{\gamma+1}{p}}\omega(S(u))^{-\f{1}{p}}.$$
It is easy to see that
$\|g_u\|_{A^p_\omega}\asymp\|h_u\|_{A^p_\omega}\asymp1$ and $g_u\to0$, $h_u\to0$  uniformly on compact subsets of $\D$ as $|u|\to1$.
Therefore,
$$\lim_{n\to\infty}\|(C_\varphi-C_\psi)g_{\varphi(z_n)}\|_{A^p_\omega}^p=0$$
and
$$\lim_{n\to\infty}\|(C_\varphi-C_\psi)h_{\varphi(z_n)}\|_{A^p_\omega}^p=0.$$
Since
$$\omega(S(z))|f(z)|^p\lesssim\|f\|_{A^p_\omega}^p,~~~\mbox{~~for all~~}~~~f\in A^p_\omega,$$
we have
$$\lim_{n\to\infty}\omega(S(z_n))\left(\left|g_{\varphi(z_n)}(\varphi(z_n))-g_{\varphi(z_n)}(\psi(z_n))\right|^p+\left|h_{\varphi(z_n)}(\varphi(z_n))-h_{\varphi(z_n)}(\psi(z_n))\right|^p\right)=0.$$
Then Lemma D yeids
\begr
\lim_{n\to\infty}\frac{\omega(S(z_n))}{\omega(S(\varphi(z_n)))}\rho(\varphi(z_n),\psi(z_n))^p=0.\nonumber
\endr

Therefore, by Lemma 1, we obtain that
\begr
\lim_{n\to\infty}\left(\frac{1-|z_n|}{1-|\varphi(z_n)|}\right)^{\b+1}\rho(\varphi(z_n),\psi(z_n))^p
=0.\nonumber
\endr
Since the two sequences $\{\frac{1-|z_n|}{1-|\varphi(z_n)|}\}$ and $\{\rho(\varphi(z_n),\psi(z_n))\}$ are both bounded.
Thus, we obtain
\begr
a=\lim_{n\to\infty}\left(\frac{1-|z_n|}{1-|\varphi(z_n)|}\right)\rho(\varphi(z_n),\psi(z_n))
=0,\nonumber
\endr
which is a desired contradiction.

Conversely, we only have to prove (10) implies that $C_\varphi-C_\psi$ is compact.
Let $\{f_k\}$ be an arbitrary bounded sequence in $A^p_\omega$ such that $f_k\to0$ uniformly on compact subsets of $\D$. It suffices to show that
$$\|(C_\varphi-C_\psi)f_k\|_{A^p_\omega}\to0,$$
as $k\to\infty$.
In order to prove this, give $0<r<1$, we put
$$E:=\{z\in\D:\rho(\varphi(z),\psi(z))<r\}~~\mbox{~~and~~}~~F:=\D\backslash E.     $$
 Then for each $k$,  
\begr
 &&\|(C_\varphi-C_\psi)f_k\|_{A^p_\omega}^p=\int_{\D}|f_k(\varphi(z))-f_k(\psi(z))|^p\omega(z) dA(z) \nonumber \\
&=& \int_E|f_k(\varphi(z))-f_k(\psi(z))|^p\omega(z) dA(z)+\int_F|f_k(\varphi(z))-f_k(\psi(z))|^p\omega(z) dA(z).  \label{di}
\endr

We first estimate the second term in the right-hand side of the equality (\ref{di}). Let $\chi_{F}$ denote the characteristic function of $F$. Since $r\chi_F\leq \rho(\varphi,\psi)$, by (\ref{cc}), we get
$$\lim_{|z|\to1}\chi_F(z)\left(\frac{1-|z|}{1-|\varphi(z)|}+\frac{1-|z|}{1-|\psi(z)|}\right)=0.$$
This, together with Lemma 5, yields
\begr
&&\int_F|f_k(\varphi(z))-f_k(\psi(z))|^p\omega(z) dA(z) \nonumber\\
&\lesssim&\int_{\D}|f_k(\varphi(z))|^p\chi_F(z)\omega(z)dA(z)
+\int_{\D}|f_k( \psi(z))|^p\chi_F(z)\omega(z)dA(z)\nonumber\\
&:=&\int_{\D}|f_k(z)|^p\nu_1(z)
+\int_{\D}|f_k(z)|^p\nu_2(z)
\to0,\nonumber
\endr
as $k\to\infty,$
where
$$\nu_1(K)=\int_{\varphi^{-1}(K)}\chi_F(z)\omega(z)dA(z)~~\mbox{~~and~~}~~
\nu_2(K)=\int_{\psi^{-1}(K)}\chi_F(z)\omega(z)dA(z),$$
for all Borel set $K\subset\D$.

Next, we estimate the first term in the right-hand side of  the equality (\ref{di}). Using Lemma 3,~Fubini's Theorem, inequality (\ref{ie}), Theorem 4  and Lemma C, we have
\begr
&&\int_E|f_k(\varphi(z))-f_k(\psi(z))|^p\omega(z) dA(z) \nonumber\\
&\lesssim&\int_{E} \rho(\varphi(z),(\psi(z)))^p\frac{\int_{\bigtriangleup(\varphi(z),r_2)}|f_k(\zeta)|^p\widetilde\omega(\zeta)dA(\zeta)}{\omega(S(\varphi(z)))}\omega(z)dA(z)
\nonumber\\
&\lesssim&r^p\int_{\D}|f_k(\zeta)|^p \frac{\int_{\varphi^{-1}(\bigtriangleup(\zeta,r_2))}\omega(z)dA(z)}{\omega(S(\zeta))}\widetilde\omega(\zeta)dA(\zeta)\nonumber\\
&\lesssim&r^p\|f_k\|_{A^p_\omega}^p\|C_\varphi\|\nonumber\\
&\lesssim&r^p.\nonumber
\endr
Letting $r\to0$, we get
$$\|(C_\varphi-C_\psi)f_k\|_{A^p_\omega}\to0.$$
The proof is complete.    $\Box$\msk

As a corollary, we obtain the following   characterization for  the operator $C_\varphi-C_\psi: A^p_\omega\to A^p_\omega$. 
The compactness of $C_\varphi-C_\psi$ on $A^p_\omega$ is independent of $p$ and $\omega$, whenever $\omega\in\mathcal{D}$.\msk

\noindent{\bf Corollary 13.}  {\it Let $0<p< \infty$ and $\omega\in\mathcal{D}$. Suppose $\varphi$ and $\psi$ are analytic self-maps of  ~$\D$. Then the  operator $C_\varphi-C_\psi: A^p_\omega\to A^p_\omega$ is compact if and only if
\begr
\lim_{|z|\to1}\left(\frac{1-|z|^2}{1-|\varphi(z)|^2}+\frac{1-|z|^2}{1-|\psi(z)|^2}\right)\rho(\varphi(z),\psi(z))=0.\nonumber
\endr}\msk

 In the rest of this section we assume that $\varphi_i:\D\to\D$ is analytic and $\varphi_i\neq\varphi_j$ if $i\neq j$. We define 
$$F_i:=\{\zeta\in\partial\D:\varphi_i \mbox{~~has a finite angular derivative at~}\zeta\}$$
and
$$\rho_{ij}(z):=\left|\frac{\varphi_i(z),\varphi_j(z)}{1-\overline{\varphi_i(z)}\varphi_j(z)}\right|.$$

The proof of the following Theorem will be quite similar to the proof of Theorem 12, with a few added
complications.\msk

\noindent{\bf Theorem 14.}  {\it Let $0<p< \infty$ and $\omega\in\mathcal{D}$.  Let
$\varphi,\varphi_1,...,\varphi_n$ be finitely many analytic self-maps of $\D$. Suppose that $C_{\varphi_1},C_{\varphi_1},\cdots,C_{\varphi_n}$ are not compact on $A^p_\omega$.
Then the  operator $C_{\varphi}-C_{\varphi_1}-\cdots-C_{\varphi_n}: A^p_\omega\to A^p_\omega$ is compact if and only if the following two conditions hold. 

(i) $F=\cup_{j=1}^nF_j$ and $F_i\cap F_j=\emptyset$ if $i\neq j$ with $i,j\geq1$;

(ii)
\begr
\lim_{z\to\zeta}\left(\frac{1-|z|^2}{1-|\varphi(z)|^2}+\frac{1-|z|^2}{1-|\varphi_j(z)|^2}\right)\rho(\varphi(z),\varphi_j(z))=0\nonumber
\endr
for all $\zeta\in F(\varphi_j)$ for $j=1,2,...,n$.
}\msk

 {\it Proof.} For the simplicity of notation, we put $T=\sum_{j=1}^nC_{\varphi_j}$.
 If $C_{\varphi}-T$ is compact on $A^p_\omega$, then by Corollary 11, $(i)$ holds.
Now, assume that $(ii)$ fails. We will derive a contradiction.

Since $(ii)$ fails, there exist $
\zeta\in F(\varphi_j)$ for some $j$ and a sequence $\{z_k\}\subset\D$ such that $z_k\to\zeta$ and
$$\lim_{k\to\infty}\rho(\varphi(z_k),\varphi_j(z_k))\left(\frac{1-|z_k|^2}{1-|\varphi(z_k)|^2}+\frac{1-|z_k|^2}{1-|\varphi_j(z_k)|^2}\right)>0.$$
 By passing to a subsequence, we may assume that
$$a_k:=\rho(\varphi(z_k),\varphi_j(z_k))\frac{1-|z_k|^2}{1-|\varphi(z_k)|^2}$$
or
$$b_k:=\rho(\varphi(z_k),\varphi_j(z_k))\frac{1-|z_k|^2}{1-|\varphi_j(z_k)|^2}$$
does not converge to zero.

Without loss of generality, we assume that $a_k$ does not converges to zero. We take 
$g_k:=g_{\varphi(z_k)}~~~~~\mbox{~~~~and~~~}h_{k}:=h_{\varphi(z_k)},$ for each $k$. 
Note that the two sequences
both $\{\rho(\varphi(z_k),\varphi_j(z_k))\}$ and $\{\frac{1-|z_k|^2}{1-|\varphi(z_k)|^2}\}$ are bounded. Thus, by  passing to anther subsequences if necessary, we may further assume
that 
$$\lim_{k\to\infty}\rho(\varphi(z_k),\varphi_j(z_k))=c_1\mbox{~~~~and~~~~}\lim_{k\to\infty}\frac{1-|z_k|^2}{1-|\varphi(z_k)|^2}=c_2,$$
for some constant $c_1,c_2>0$ with $c_1\leq1$.

Also, note that $\zeta\notin F(\varphi_i)$ for $i\neq j$. By the Julia-Caratheodory Theorem, we have
$$\lim_{k\to\infty}\frac{1-|z_k|}{1-|\varphi_i(z_k)|}=0,~~~~i\neq j,$$

\begr
\lim_{k\to\infty}\omega(S(z_k))|g_k(\varphi_i(z_k))|^p&=&\lim_{k\to\infty}\frac{\omega(S(z_k))}{\omega(S(\varphi_i(z_k)))}\left|\frac{1-|\varphi(z_k)|^2}{1-\overline{\varphi(z_k)}\varphi_i(z_k)}\right|^{\gamma+1}
\nonumber\\
&\lesssim&\lim_{k\to\infty}\left(\frac{1-|z_k|}{1-|\varphi_i(z_k)|}\right)^{\a+\gamma+2}=0.\nonumber
\endr

\begr
\lim_{k\to\infty}\omega(S(z_k))|h_k(\varphi_i(z_k))|^p&=&\lim_{k\to\infty}\frac{\omega(S(z_k))}{\omega(S(\varphi_i(z_k)))}\left|\frac{1-|\varphi(z_k)|^2}{1-t_N\overline{\varphi(z_k)}\varphi_i(z_k)}\right|^{\gamma+1}
\nonumber\\
&\lesssim&\lim_{k\to\infty}\left(\frac{1-|z_k|}{1-|\varphi_i(z_k)|}\right)^{\a+1}\left(\frac{1-|z_k|}{1-t_N|\varphi_i(z_k)|}\right)^{\gamma+1}\nonumber\\
&\lesssim&\lim_{k\to\infty}\left(\frac{1-|z_k|}{1-|\varphi_i(z_k)|}\right)^{\a+\gamma+2}=0.\nonumber
\endr

The same argument as in the proof of Theorem 12 yields
$$\lim_{k\to\infty}\omega(S(z_k))\left(|g_k(\varphi(z_k))-(Tg_k)(z_k)|^p+|h_k(\varphi(z_k))-(Th_k)(z_k)|^p\right)=0.$$
Thus, the same argument as in the proof of Theorem 12 yields
\begr
\lim_{k\to\infty}\left(\frac{1-|z_k|}{1-|\varphi(z_k)|}\right)\rho(\varphi(z_k),\varphi_j(z_k))
=0, \nonumber
\endr
which is a desired contradiction.

Assume next that both $(i)$ and $(ii)$ hold. We will prove that $C_{\varphi}-T$ is compact. The proof will be quite similar to the proof of Theorem 12. 
~Define $$D_i :=\Big \{z\in\D : \frac{1-|z|^2}{1-|\varphi_i(z)|^2}\geq \frac{1-|z|^2}{1-|\varphi_j(z)|^2},\mbox~{~for~all~} j\neq i\Big\}$$~
for $i = 1, . . . , N$. Fix $0<r<1$ and define
$$E_i:=\{z\in D_i: \rho(\varphi(z),\varphi_i(z))<r\}~ \mbox{~~~and~~~}~ E_i^\prime:=D_i\backslash E_i.$$
By the proof of \cite[Theorem 5]{Mo}, we get
\begin{equation}
\lim_{|z|\to1}\chi_{E_i^\prime}(z)\left(\frac{1-|z|}{1-|\varphi(z)|}+\frac{1-|z|}{1-|\varphi_j(z)|}\right)=0,~~\mbox{~for all~~}i,j,  \label{mot}
\end{equation}
and
\begin{equation}
\lim_{|z|\to1}\chi_{E_i}(z)\frac{1-|z|}{1-|\varphi_j(z)|}=0,~~\mbox{~whenever~}i\neq j.\label{mot2}
\end{equation}

Now, let $\{f_n\}$ be a bounded sequence in $A^p_\omega$ such that $f_k\to0$ uniformly on compact subset of $\D$. Since $\D=\cup_{i=1}^nD_i$, we have
$$\|(C_\varphi-T)f_k\|_{A^p_\omega}^p=\int_{D}|f_k\circ\varphi-\sum_{i=1}^nf_k\circ\varphi_i|^p\omega dA
\leq\sum_{i=1}^n\int_{E_i}+\sum_{i=1}^n\int_{E_i^\prime}.$$
Note, as in the proof of Theorem 12, that the second sum of the above tends to $0$ as $k\to\infty$, by equality (\ref{mot})~and Lemma 5.
For the $i$-th term of the first sum, we have
$$\int_{E_i}\lesssim\int_{E_i}|f_k\circ\varphi-f_k\circ\varphi_i|^p\omega dA+\sum_{j\neq i}\int_{E_i}|f_k\circ\varphi_j|^p\omega dA.$$
Note from equality (\ref{mot2}) and Lemma 5 that the second term of the above tends to 0 as $k\to\infty$.
Finally, from the proof of Theorem 12 we see that the  first term of the above is dominated by $r^p$.
So, we conclude that
$$\limsup_{k\to\infty}\|(C_\varphi-T)f_k\|_{A^p_\omega}^p\lesssim r^p.$$
Letting $r\to0$, we obtain
$\limsup_{k\to\infty}\|(C_\varphi-T)f_k\|_{A^p_\omega}^p=0.$
The proof is complete.$\Box$\msk

Theorem 14 and Corollary 9 immediately yield the following characterization for a composition operator to be equal module compact operators to a linear combination of composition operators.\msk

\noindent{\bf Theorem 15.}  {\it Let~$0<p< \infty$ and $\omega\in\mathcal{D}$. Let
$\varphi,\varphi_1,...,\varphi_n$ be finitely many analytic self-maps of $\D$. Suppose that $C_{\varphi}$,~$C_{\varphi_1},\cdots,C_{\varphi_n}$ are not compact on $A^p_\omega$. Let $\lambda_1,\cdots,\lambda_n\in \C\backslash\{0\}$.    Then the  operator $C_{\varphi}-\sum_{j=1}^n\lambda_jC_{\varphi_j}: A^p_\omega\to A^p_\omega$ is compact if and only if the following three conditions holds:

(1) $\lambda_1=...=\lambda_n=1;$

(2) $F=\cup_{j=1}^nF_j$ and $F_i\cap F_j=\emptyset$ if $i\neq j$ with $i,j\geq1$;

(3)
\begr
\lim_{z\to\zeta}\left(\frac{1-|z|^2}{1-|\varphi(z)|^2}+\frac{1-|z|^2}{1-|\varphi_j(z)|^2}\right)\rho(\varphi(z),\varphi_j(z))=0\nonumber
\endr
for all $\zeta\in F_j$ for $j=1,2,...,n$.
}\msk


\begin{thebibliography}{99999}

\bibitem{Be} E. Berkson, Composition operators isolated in the uniform operator topology,  {\it Proc. Amer. Math. Soc.}  {\bf 81(2)} (1981), 230--232.

\bibitem{CKP2} B. Choe, H. Koo and I. Park, Compact differences of composition operators over polydisks, {\it
Integr. Equ. Oper. Theory} {\bf 73} (2012) 57--91.

\bibitem{CKP1} B. Choe, H. Koo and I. Park, Compact differences of composition operators on the Bergman spaces over the ball, {\it Potential Anal.}
 {\bf 40(1)} (2014), 81--102.



\bibitem{cm} C. Cowen and B. Maccluer, {\it Composition Operators on Spaces of Analytic Functions}, CRC Press, Boca Raton, FL, 1995.


\bibitem{DW} P. Duren and R. Weir, The pseudo-hyperbolic metric and Bergman spaces in the ball, {\it Trans. Amer. Math. Soc.} {\bf 359(1)} (2007), 63--76.

\bibitem{G} T. Goebeler, Composition operators acting between Hardy spaces, {\it Integr. Equat. Oper. Th.} {\bf 41(4)} (2001), 389--395.

\bibitem{KM} T. Kriete and J. Moorhouse, Linear relations in the calkin algebra for composition operators, {\it Trans. Amer. Math. Soc.} {\bf 359(6)} (2007), 2915--2944.

\bibitem{KW} H. Koo and M. Wang, Joint Carleson measure and the difference of composition operators on $A^p_\a(\B_n)$,
 {\it J. Math. Anal. Appl.} {\bf 419(2)} (2014), 1119--1142.


\bibitem{LS} M. Lindstrom and E. Saukko, Essential norm of weighted composition operators and difference of composition operators
between standard weighted Bergman spaces,  {\it Complex Anal. Oper. Theory} {\bf 9(6)} (2015), 1411--1432.


\bibitem{Mo} J. Moorhouse, Compact differences of composition operators, {\it J. Funct. Anal.}  {\bf 219} (2005), 70--92.

\bibitem{NS} P. Nieminen and E. Saksman, On compactness of the difference of composition operators, {\it J. Math. Anal. Appl.} {\bf  298(2)}(2004),   501--522.

\bibitem{P} J. Pel{\'a}ez, Small weighted Bergman spaces, {\it In: Proceedings of the summer school in ``complex and
harmonic analysis", and related topics}, (2016), 29--98.

\bibitem{PR1} J. Pel{\'a}ez and  J. R{\"a}tty{\"a}, Weighted Bergman spaces induced by rapidly increasing weights, {\it Mem. Am. Math. Soc.}, {\bf 227(1066)}(2014).

\bibitem{PR2} J. Pel{\'a}ez and J. R{\"a}tty{\"a}, Embedding theorems for Bergman spaces via harmonic analysis, {\it Math. Ann.} {\bf 362(1-2)} (2015), 205--239.


\bibitem{PR3} J. Pel{\'a}ez and R. Jouni, Trace class criteria for Toeplitz and composition operators on small Bergman spaces, {\it Adv. Math.} {\bf 293} (2016), 606--643.

\bibitem{PR4} J. Pel{\'a}ez and J. R{\"a}tty{\"a}, Two weight inequality for Bergman projection, {\it J. Math. Pures Appl.} {\bf 105(1)} (2016), 102--130.


\bibitem{PR5} J. Pel{\'a}ez and J. R{\"a}tty{\"a}, Hankel operators induced by radial Bekoll{\'e}-Bonami weights on Bergman spaces, {\it Math. Z.} (2019) 1--28.


\bibitem{PRS} J. Pel{\'a}ez, J. R{\"a}tty{\"a} and K. Sierra, Embedding Bergman spaces into tent spaces, {\it Math. Z.} {\bf 281(3-4)} (2015), 1215--1237.

\bibitem{S1} E. Saukko, Difference of composition operators between standard weighted Bergman spaces, {\it J. Math. Anal. Appl.} {\bf  381(2)} (2011),   789--798.

\bibitem{S2} E. Saukko,  An application of atomic decomposition in Bergman spaces to the study of differences of composition operators, {\it J. Funct. Anal.} {\bf 262(9)} (2012), 3872--3890.

\bibitem{SLD} Y. Shi, S. Li and J. Du, Difference of composition operators between weighted Bergman spaces on the unit ball, arXiv::1903.00651

\bibitem{Sh1} J. Shapiro, {\it Composition Operators and Classical Function Theory}, Springer Science,  Business Media, 2012.

\bibitem{SS} J. Shapiro and C. Sundberg, Isolation amongst the composition operators, {\it Pacific J. Math.}  {\bf 145} (1990), 117--152.

 \bibitem{SL} Y. Shi and S. Li, Difference of composition operators between different Hardy spaces, {\it J. Math. Anal. Appl.} {\bf 467(1)} (2018), 1--14.

\bibitem{Zhu}  K.~Zhu,  {\it Operator Theory in Function Spaces},  American Mathematical Society, Providence, RI, 2007.



\end{thebibliography}
\end{document}